\documentclass[12pt]{article}

\newcommand  {\Ebar} {{\mbox{\rm$\mbox{I}\!\mbox{E}$}}}
\newcommand  {\Rbar} {{\mbox{\rm$\mbox{I}\!\mbox{R}$}}}

\newcommand  {\Nbar} {{\mbox{\rm$\mbox{I}\!\mbox{N}$}}}

\newsavebox{\zzzbar}
\sbox{\zzzbar}
  {\setlength{\unitlength}{0.9em}
  \begin{picture}(0.6,0.7)
  \thinlines
  \put(0,0){\line(1,0){0.6}}
  \put(0,0.75){\line(1,0){0.575}}
  \multiput(0,0)(0.0125,0.025){30}{\rule{0.3pt}{0.3pt}}
  \multiput(0.2,0)(0.0125,0.025){30}{\rule{0.3pt}{0.3pt}}
  \put(0,0.75){\line(0,-1){0.15}}
  \put(0.015,0.75){\line(0,-1){0.1}}
  \put(0.03,0.75){\line(0,-1){0.075}}
  \put(0.045,0.75){\line(0,-1){0.05}}
  \put(0.05,0.75){\line(0,-1){0.025}}
  \put(0.6,0){\line(0,1){0.15}}
  \put(0.585,0){\line(0,1){0.1}}
  \put(0.57,0){\line(0,1){0.075}}
  \put(0.555,0){\line(0,1){0.05}}
  \put(0.55,0){\line(0,1){0.025}}
  \end{picture}}
\newcommand{\Zbar}{\mathord{\!{\usebox{\zzzbar}}}}

\newsavebox{\uuunit}
\sbox{\uuunit}
    {\setlength{\unitlength}{0.825em}
     \begin{picture}(0.6,0.7)
        \thinlines
        \put(0,0){\line(1,0){0.5}}
        \put(0.15,0){\line(0,1){0.7}}
        \put(0.35,0){\line(0,1){0.8}}
       \multiput(0.3,0.8)(-0.04,-0.02){12}{\rule{0.5pt}{0.5pt}}
     \end {picture}}

\newcommand{\QED}{{\hspace*{\fill}\rule{2mm}{2mm}\linebreak}}

\newtheorem{lemma}{Lemma}[section]
\newtheorem{proposition}{Proposition}[section]
\newtheorem{theorem}{Theorem}[section]

\newtheorem{corollary}{Corollary}[section]
\newcommand{\Z}{\Zbar}
\newcommand{\R}{\Rbar}
\newcommand{\N}{\Nbar}

\newcommand{\E}{\Ebar}

\begin{document}
\setlength{\textheight}{21cm}
\title{On the thermodynamic limit for a one-dimensional sandpile process}
\author{C. Maes\thanks{Onderzoeksleider FWO, Flanders. \
                        Email: christian.maes@fys.kuleuven.ac.be }  \\
    F. Redig\thanks{Post-doctoraal onderzoeker FWO, Flanders. \
    Email: frank.redig@fys.kuleuven.ac.be} \\
    E. Saada\thanks{CNRS, UPRESA 6085, Rouen, France.\
    Email: Ellen.Saada@univ-rouen.fr}\\
    A. Van Moffaert\thanks{Aspirant FWO, Flanders. \
                       Email: annelies.vanmoffaert@fys.kuleuven.ac.be }
                 \thanks{Address :
    Celestijnenlaan 200D, B-3001 Leuven, Belgium } }
\date{ }
\maketitle

\begin{abstract}
Considering the standard abelian sandpile model in one dimension,
we construct an infinite volume Markov process corresponding  to
its thermodynamic (infinite volume) limit.  The main difficulty we
overcome is the strong non-locality of the dynamics.  However,
using similar ideas as in recent extensions of the standard Gibbs
formalism for lattice spin systems, we can identify a set of `good'
configurations on which the dynamics is effectively local.  We
prove that every configuration converges in a finite time to the
unique invariant measure.
\end{abstract}
\vspace{3mm}
\begin{flushleft} 


\end{flushleft}

\section{Introduction.}
The standard abelian sandpile model was introduced in 1988,
\cite{Bak}.   This model lives on a finite graph (e.g. a square on
$\Z^2$) in which to each vertex is associated a height-variable
(representing the height of a sand column at that site or, the
averaged difference in sandlevels with neighboring sites). A vertex
is picked at random and its height is increased by one. If the
height variable exceeds some critical value, then the vertex
becomes ``unstable'' and ``topples'', i.e. gives an equal portion
of its sand to each of its neighbours (adjacent vertices) which in
their turn, can become ``unstable'' and ``topple'', and so on,
until every vertex has again a subcritical height-value. Some
vertices are connected to a sink (the boundary) where sand
disappears. An unstable vertex thus creates an ``avalanche'' which
can cause the toppling of many vertices around it. The number of
vertices involved in or affected by one avalanche can be
arbitrarily large (dependent on the configuration). 

Over the last
decade, sandpile processes and various versions thereof have been
intensively studied as simple lattice models exhibiting the
phenomenon of self-organized criticality.  While it is still a
subject of debate what is the precise nature of the phenomenon and
under exactly what conditions it obtains, sandpiles have become
widely appreciated as simple threshold dynamics showing scale
invariance of various response functions with (critical) exponents
of an astounding universality. Not surprisingly, a vast amount of
computer simulations have been performed involving sandpile
processes and their acclaimed self-organized criticality. Yet, both
controlled laboratory or real nature experiments and a rigorous
mathematical analysis have remained rather limited.  We refer to
\cite{Strauven} and the references therein for a number of
experiments validating some of the universal aspects of the
sandpile paradigm. In this paper we do not address questions of
self-organized criticality nor do we refer to the most interesting
physical realizations of sandpile-like models.  Rather we turn to
the simplest mathematical questions one can ask from the point of
view of the theory of interacting particle systems for the simplest
one-dimensional abelian sandpile process.

 For the abelian sandpiles the first mathematical results
 have been obtained by D. Dhar, see \cite{Dhar1, Dhar2, Dhar3} and the review
 paper \cite{Dhar4}.  As further discussed
 also in \cite{Speer},  the main question there was
the characterization of 
the invariant measure and the recurrent configurations
of the finite Markov chain defined in the sandpile process.
In one dimension, the answers are rather simple.  
A detailed study can be found
in \cite{ruelle}.  Our questions concern the so called thermodynamic limit,
i.e. how to construct a {\it bona fide} infinite volume limit 
dynamics and how to characterize
its main properties. 
In \cite{Priezz}, Priezzhev was able to compute exactly the
single-site marginal of the invariant measure of the sandpile
process on a square in the limit as the square tends to cover the
whole plane.  Various other mathematical results have been obtained
but, for all we know, the problem of the thermodynamic limit has
remained widely open so far.  In the present paper we give a
complete solution to this problem for the standard abelian sandpile
process in one dimension. 

It should be realized that this
one-dimensional model has as such rather limited interest from the
physics point of view.  Yet, some of the mathematical questions
pertaining to the thermodynamic limit remain very non-trivial and,
arguably have answers which are expected to be not so different in
higher dimensions.  It is therefore unavoidable to start in one
dimension and clarify the situation there.  The main mathematical
difficulty remains unchanged and has to do with the strong
non-locality of the dynamics.  It is built in in the definition of
the sandpile process and it seems to be essential for some of its
most interesting aspects.  The effect of adding a sand particle at
one site can affect the whole system.
 This implies that the standard machinery to
construct a Feller process (Hille-Yoshida theorem, martingale
problem) does not work.  Similar problems were encountered in the
study  of long range exclusion processes, see \cite{Ligg1},
\cite{Herve}, \cite{Carlson1}.  Here the situation is even more
exciting as will be seen from the convergence to the stationary
measure in a finite time and from the relation between the process
and its formal generator.  This relation is only `pointwise' for a
certain class of `decent' configurations (in particular the process
does not have a generator in the standard sense).  One finds
analogies with some recent work on weakly and almost Gibbsian
measures where the interaction energy is only absolutely summable
on a full measure set of configurations but not uniformly, see
\cite{MRV1,MRV2}.

Our paper is organized as follows. In Section 2 we introduce the
model and the formal generator of the sandpile process. In Section
3 we prove how to construct a Markov semigroup, and thus a Markov
process which corresponds to this formal generator. The main
technical tool in this construction is monotonicity. In Section 4
we show the precise relation between the semigroup and the formal
generator. In Section 5 we study the invariant measure of the
process and show that every initial configuration converges in a
finite time to this invariant measure.

\section{Definitions.}
\subsection{Configurations.}
The state space of the one-dimensional 
sandpile process is $\Omega=\{1,2\}^{\Z}$. Elements of $\Omega$ are called
configurations and are denoted by $\eta,\xi,\zeta$.
The state space is a compact metric space in the product
topology. 
The value $\eta (x) = 2$ is called the {\it threshold value}. Sites for
which $\eta (x)=2$ are called {\it critical sites} in the configuration $\eta$.
The set of all finite subsets of $\Z$ is denoted by $\mathcal{S}$. 
We define $\Omega_f$ to be the set of configurations with a
finite number of critical sites, i.e.
\begin{equation}
\Omega_f := \{ \eta \in \Omega: \eta^{-1} (\{2\})\in \mathcal{S}\}.
\end{equation}
There is a one-to-one correspondence between $\mathcal{S}$ and $\Omega_f$
via $A\rightarrow\eta_A$, where $\eta_A$ denotes the configuration whose
set of critical sites is $A$, i.e. $\eta_A^{-1} (\{2\})=A$.
The configurations in $\Omega_f$ will play an important role in the
construction of the process (cf. below). A function $f:\Omega\rightarrow \R$
is called local if it depends only on a finite number of coordinates.
Every local function is continuous and every continuous function is a uniform
limit of local functions.

\subsection{Toppling transformation $T_i$.}
Given a configuration $\eta\in \Omega_f$ and a lattice site $i\in \Z$,
$T_i\eta$ represents the configuration obtained by adding one sand particle
at site $i$ and letting the system {\it ``topple"} until a (stable) 
configuration
(i.e. an element of $\Omega_f$) is obtained. 
This toppling is described in words 
as follows:
If $\eta(i)=2$ and one sand particle is added at site $i$, then from that site
two particles will be removed, one to the left neighbour $i-1$ and one to
the right neighbour $i+1$. If say 
$\eta (i+1)$ was already $2$, then site $i+1$ will
also ``topple", giving one particle to $i$, and one to $i+2$. This
goes on until a configuration in which 
$\eta(j)\leq 2$ for all $j\in \Z$ is obtained. 
This succession of topplings at different
sites caused by adding one particle to a critical site is called an
{\it avalanche}, and thus $T_i\eta$ represents the configuration after the
avalanche. If $\eta (i) =1$, then $T_i\eta$ differs from $\eta$
at site $i$ only ($T_i\eta (i) = 2$), but if $\eta (i)=2$, then 
$T_i\eta$ differs from
$\eta$ at three lattice sites which depend on $i$ and $\eta$ in a non-local
way.
This corresponds exactly to the dynamics of the standard abelian sandpile
process in a finite lattice region.

More precisely we introduce
\begin{equation}
k^+ (i,\eta):= \inf \{j\geq 0:\eta (i+j) =1 \},
\end{equation}
and
\begin{equation}
k^- (i,\eta) := \inf \{ j>0: \eta (i-j) =1 \},
\end{equation}
where $\inf\emptyset := +\infty$.
We denote by $e_i$ the mapping from $\Z$ to $\{0,1\}$ which is one
at site $i$ and zero at all other site. We 
distinguish five cases:
\begin{enumerate}
\item $k^+ (i,\eta)=0$, i.e. $\eta(i)=1$, then
\begin{equation}
T_i\eta = \eta + e_i.
\end{equation}
\item $k^+(i,\eta ) >0$, $k^+ (i,\eta) \vee k^- (i,\eta ) <\infty$, then
\begin{equation}
T_i\eta = \eta + e_{i+ k^+(i,\eta)}+e_{i-k^-(i,\eta)}
-e_{i+k^+(i,\eta)-k^-(i,\eta )}
\end{equation}
In what follows we abbreviate
\begin{equation}
i^\eta:= i + k^+(i,\eta )-k^- (i,\eta ).
\end{equation}
\item $k^+ (i,\eta ) =\infty, k^- (i,\eta) <\infty$, then
\begin{equation}
T_i\eta = \eta + e_{i-k^-(i,\eta)}
\end{equation}
\item $k^+ (i,\eta)<\infty, k^- (i,\eta ) =\infty$, then
\begin{equation}
T_i\eta = \eta + e_{i+k^+ (i,\eta )}
\end{equation}
\item $k^+(i,\eta ) = k^- (i,\eta ) = \infty$ (i.e. $\eta \equiv 2$), 
then
\begin{equation}\label{conv}
T_i\eta = \eta.
\end{equation}
\end{enumerate}

Cases 3,4,5 above are just conventions to define $T_i\eta$ for
all $\eta\in\Omega$, and to write down the formal generator,
but they are not relevant in the construction of the sandpile
process.
That these choices are `correct' will be seen later once it is shown that the
configuration $\eta \equiv 2$
is invariant for the constructed process (cf. Proposition \ref{deltatwee}).

The explicit and simple expression of $T_i\eta$ is typical for
the model in dimension one. In higher dimensions $T_i \eta$ will become
a much more complicated function of $\eta$.
Note that the transformations $T_i$ and $T_j$ commute, i.e.
\begin{equation}
T_i(T_j(\eta)) = T_j (T_i(\eta)), \ \forall i,j\in \Z,\ \forall \eta\in \Omega.
\end{equation}
For this reason the sandpile model has been called {\it abelian}
in the physics literature, cf. \cite{Dhar1}.

\subsection{Monotonicity.} 
In the construction of the sandpile process, monotonicity will play a crucial
role, just as in the construction of the long-range exclusion process,
cf. \cite{Ligg1}, \cite{Herve}. For more details on monotonicity, we
refer to \cite{Ligg2}, Chapter 2, Section 2.

For $\eta,\xi\in \Omega$ we define $\eta\leq\xi$ if $\eta (x) \leq\xi (x)$
for all $x\in\Z$. A function $f:\Omega\rightarrow\R$ is called {\it monotone}
if $\eta\leq\xi$ implies $f(\eta ) \leq f(\xi )$, for all 
$\eta,\xi \in\Omega$. We denote by $\mathcal{M}$ the class of all monotone
bounded Borel measurable functions.
A Markov process $\{ \eta_t :t\geq 0 \}$ on $\Omega$ with Markov semigroup
$\{ S(t):t\geq 0\}$ is called monotone if $f\in\mathcal{M}$ implies
$S(t)f \in \mathcal{M}$. Usually one proves monotonicity of a Markov process
by explicitly constructing a coupling of path space measures such that for
$\eta\leq\xi$ one has $P^{\eta,\xi} (\eta_t\leq\xi_t,\  \forall t\geq 0) = 1$.

\subsection{Formal generator.}
In sandpile processes, the avalanches are immediate. Once the vertex $i$ is
chosen, the transformation $\eta \rightarrow T_i \eta$
happens instantaneously. This is 
what one calls {\it infinite 
separation of time scales}
in the physics literature, and it has been argued as 
essential to have a so-called 
self-organized critical state.

The formal generator of the sandpile process has therefore to be written as
\begin{equation}\label{watte}
Lf(\eta ) = \sum_{i\in \Z} [f(T_i\eta)- f(\eta )].
\end{equation}
Since $L$ is of course not the generator of a Feller semigroup, 
we shall have
to give a precise meaning to this description of the dynamics, and to clarify 
how to ``associate" a process to $L$.

We can rewrite the formal generator
\begin{equation}
Lf(\eta ) = L_b f(\eta ) + L_a f(\eta ),
\end{equation}
where
\begin{equation}
L_b f(\eta ) = \sum_{i\in\Z} \chi (\eta(i) = 1 ) [f(\eta+\delta_i ) -f(\eta)]
\mbox{ is the ``birth part"},
\end{equation}
\begin{equation}
L_a f(\eta ) = \sum_{i\in\Z} \chi (\eta (i) = 2) [f(T_i\eta) - f(\eta)]
\mbox{ is the ``avalanche part"},
\end{equation}
and $\chi$ denotes the indicator function. 
The operator 
$L_b$ corresponds to a pure birth process (independent births with
rate one at each lattice site) and thus is a well-defined generator of
a Feller semigroup. The avalanche 
part however gives technical problems related
to the nonlocality of the transformation $T_i$. One cannot 
run the standard Hille-Yoshida program, which would give a Feller semigroup
associated to $L_b+L_a$. Indeed, the process we shall associate to $L$
 will be constructed by monotonicity
and turn out not to be Feller. 

{\bf Remark.} The processes we consider here run in continuous time while it is
a standard procedure to describe the finite volume sandpile process
in discrete time. However, since both updating mechanisms are sequential, the
difference is unimportant.

\section{Construction of the sandpile process.}
We proceed in three steps:
\begin{enumerate}
\item Definition of a process associated to $L_a$ (only avalanches) on
$\Omega_f$. Construction of a coupling showing that this process is monotone.
\item Definition of a process associated to $L_a$ + births in the finite
interval $[-n,n]$ on $\Omega_f$. Construction of a coupling showing that
this process is still monotone. This will give us a semigroup defined on
$f\in \mathcal{M}$ by 
$S_n(t) f(\eta) = \lim_{\eta'\in\Omega_f,\ \eta'\uparrow\eta} S_n(t) f(\eta')$,
where $\lim_{\eta'\uparrow\eta}$ denotes the limit along sequences
$\{\eta'_n, n\in\N\}\subset\Omega_f$ such that ${\eta'}_n\leq{\eta'}_{n+1}$.
\item Monotonicity of the semigroups $S_n(t)$ in $n$, i.e., 
for all $f\in \mathcal{M}$, $ \eta \in \Omega$, $ t\geq 0$,
$S_{n+1} (t) f(\eta ) \geq S_n (t) f (\eta )$.
\end{enumerate}
We shall finally define the sandpile process as the process associated
to the semigroup
\begin{equation}
S(t)f(\eta) = \lim_{n\uparrow\infty}\lim_{\eta'\in \Omega_f, \eta'\uparrow\eta}
S_n(t)f(\eta').
\end{equation}

\subsection{Step 1: Avalanche part of the generator.}
For $\eta\in \Omega_f$, $A(\eta ) := \eta^{-1} (\{2\})\in\mathcal{S}$ and
we can construct the process with generator $L_a$ from the non-exploding 
continuous time Markov
chain $\{ A_t:t\geq 0\}$ on $\mathcal{S}$ with generator
\begin{equation}
\hat{L_a} f(A) 
:= \sum_{i\in A} \left[
f(A\setminus \{i^{\eta} \}
\cup \{i+k^+(i,\eta)  \} \cup \{i-k^-(i,\eta )\})
-f(A) \right],
\end{equation}
for $f:\mathcal{S} \rightarrow \R$ bounded.
Indeed, 
$\E^A |A_t| = |A| e^t$ (where $|A|$ denotes the cardinality of $A$)
so that $A_t\in \mathcal{S}$, for all $t\geq 0$,
$P^A$-a.s.
The process with generator $L_a$ (avalanche part) is thus defined on $\Omega_f$
by $\eta_t := \eta_{A_t}$, where $\{ A_t:t\geq 0 \}$ starts from $A_0 = \eta^{-1} (\{ 2\} )$. We 
show that this process is monotone via the coupling of $\eta_t,\xi_t$ ( for 
$\eta, \xi \in \Omega_f$ with $\xi\leq\eta$) described as follows. When a 
sand particle is added in
$\eta$ at site $i\in \Z$ where $\eta(i)=2$, an avalanche
creates a $1$ at site $i^\eta=i+k^+(i,\eta)-k^-(i,\eta)$. Define now 
$\phi (i,\eta,\xi)$ as the unique site on $\Z$ such that 
a sand particle added at that site in 
configuration $\xi$ would
create a $1$ in $\xi$ at the same site $i^\eta$. More precisely,
if $\xi (i^\eta ) = 2$, 
\begin{equation}
\phi(i,\eta,\xi )= i^\eta + k^+ (i^\eta,\xi ) - k^- (i^\eta, \xi ),
\end{equation}
and if $\xi (i^\eta ) =1$, we define $\phi (i,\eta,\xi ) = \infty$ 
and
$T_\infty \xi =\xi$. We now write the coupling generator 
on $\Omega_f\times\Omega_f$
for initial configurations $(\eta,\xi)$ with $\eta\geq\xi$ as
\begin{equation}
L_a^c f(\eta,\xi ) = \sum_{i\in\Z} \chi (\eta(i) = 2)
[f(T_i\eta, T_{\phi (i,\eta,\xi )}\xi ) -f(\eta,\xi )].
\end{equation}

\begin{proposition}
\begin{enumerate}
\item On $\Omega_f\times\Omega_f$ the generator $L^c_a$ defines a Markov 
process. 
For $\eta\geq \xi$, the path space measure $P^{\eta,\xi}$ is a coupling
of $P^\eta$ and $P^\xi$ (the path space measures of the processes with generator
$L_a$ starting respectively from $\eta$ and $\xi$).
\item For $\eta\geq\xi$ we have 
$P^{\eta,\xi} (\eta_t \geq \xi_t,\ \forall t\geq 0 )=1$.
\item The process with generator $L_a$ is monotone on $\Omega_f$. For the
semigroup $S_a (t):=\exp(tL_a)$ 
associated to $L_a$, $f\in \mathcal{M}$ implies
$S_a (t) f \in \mathcal{M}$.
\end{enumerate}
\end{proposition}
{\bf Proof:} 
To check that the process with generator
$L^c_a$ has the
right second marginal for $\eta\geq\xi$,  
remark that
for every $j\in\Z$ such that $\xi (j)=2$, there exists a unique
$i\in\Z$ such that $\phi (i,\eta,\xi )=j$
($L_a^c$ does not define a coupling of $P^\eta$ and $P^\xi$ for
general couples $(\eta,\xi)$ but only for the ordered couples, i.e. 
$(\eta,\xi)$ with $\eta\geq\xi$ and $\eta, \xi\in \Omega_f$). For points
2 and 3, note
that by construction any transition in the coupled process preserves the
order.\QED

\subsection{Step 2: Avalanche + birth in a finite interval.}
We now include the birth part of the formal generator $L$, only on the
interval $[-n,n]$. We have to prove 
that the associated process on $\Omega_f$ is
still monotone, since births and
avalanches do not commute.
For this, the generator
\begin{eqnarray}
L^{(n)} f (\eta ) &=& \sum_{i\in \Z} \chi (\eta (i) = 2)
[f(T_i\eta ) - f(\eta )]\nonumber\\ 
&+& \sum_{i=-n}^n \chi (\eta (i) = 1) [f(T_i\eta ) - f(\eta )]
\end{eqnarray}
defines a continuous time Markov chain
on $\mathcal{S}$, and thus a process on $\Omega_f$. We now couple two versions
of this process starting from $\eta,\xi\in \Omega_f$ such that $\eta\geq\xi$,
with the generator
\begin{eqnarray}
L^{c,(n)} f (\eta, \xi ) &=&
\sum_i \chi (\eta (i) = 2) \chi (\xi (i^\eta) = 2)
[f(T_i\eta,T_{\phi(i,\eta,\xi )}\xi)-f(\eta,\xi)]
\nonumber\\
&+& \sum_{i\in [-n,n]} \chi (\eta (i) = 2 ) \chi (\xi(i) = 1)
[f(\eta,T_i\xi)-f(\eta,\xi)]\nonumber\\
&+&\sum_{i\in[-n,n]}\chi(\eta (i) = 1 ) \chi (\xi (i) =1)
[f(T_i\eta,T_i\xi) - f(\eta,\xi ) ]\nonumber\\
&+& \sum_i \chi (\eta (i) = 2 ) \chi (\xi (i^\eta) =1)
[f(T_i\eta,\xi ) - f(\eta,\xi)].
\end{eqnarray}
This coupling is described as follows: To each site $i\in \Z$ we
associate two independent Poisson clocks, C and C'. C indicates
the ``avalanche event-times" and C' the ``birth event-times".
When C rings at $i$ and $\eta(i)=2$, we simultaneously add a sand particle
at $i$ in $\eta$ and at $\phi(i,\eta,\xi)$ in $\xi$ (if
$\phi (i,\eta,\xi ) = \infty$, we have no transition in $\xi$).
When C' rings at $i$ and $\eta (i ) =1$, then,
since $\eta\geq\xi$, $\xi (i)=1$ and we add 
a sand particle in both $\xi$ and $\eta$
simultaneously. If however $\eta (i) = 2$ and $\xi (i) = 1$, then we only
add in $\xi$ at site $i$. Like in Section 3.1, this construction implies

\begin{proposition}
\begin{enumerate}
\item The generator $L^{c,(n)}$ defines a Markov process on $\Omega_f\times\Omega_f$.
For $\eta\geq\xi$ its path space measure $P_n^{(\eta,\xi)}$ is a coupling
of $P^\eta_n$ and $P^\xi_n$, the path space measures of the processes with
generator $L^{(n)}$ starting respectively from $\eta$ and $\xi$.
\item If $\eta\geq\xi$, then 
$P_n^{(\eta,\xi)} (\eta_t\geq\xi_t,\ \forall t\geq 0)
=1$.
\item The process with generator $L^{(n)}$ is monotone on $\Omega_f$. For
the associated semigroup $S_n(t):=\exp(tL^{(n)})$, $f\in\mathcal{M}$ implies
$S_n(t) f\in \mathcal{M}$.
\end{enumerate}
\end{proposition}

\subsection{Step 3: Avalanche + births in different intervals.}
This final step is the easiest: By allowing births on a
larger interval, we can only obtain larger configurations.
We couple the processes with generators $L^{(n+1)}$ and $L^{(n)}$ through the
following generator, where $\eta, \xi \in \Omega_f$, $\eta\geq\xi$.
\begin{eqnarray}
L^{c,n+1,n} f (\eta, \xi ) &=&
\sum_i \chi (\eta (i) = 2 ) \chi (\xi (i^\eta ) = 2))
[f(T_i\eta,T_{\phi(i,\eta,\xi)})-f(\eta,\xi)]\nonumber\\
&+&\sum_{i\in [-n,n]} \chi (\eta (i)=2) \chi (\xi (i) = 1 )
[f(\eta,T_i\xi ) - f(\eta,\xi ) ]\nonumber\\
&+& \sum_{i\in [-n,n]} \chi (\eta (i) = 1 ) \chi (\xi (i) = 1)
[f(T_i\eta,T_i\xi ) - f(\eta,\xi)]\nonumber\\
&+& \sum_{i\in \{-n-1,n+1\}} \chi (\eta (i) = 1)\chi (\xi (i) = 1)
[f(T_i\eta,\xi ) - f(\eta, \xi )]\nonumber\\
&+& \sum_i \chi (\eta (i) = 2) \chi (\xi (i^\eta) =1 ) 
[f(T_i\eta,\xi ) - f (\eta, \xi ) ].
\end{eqnarray}

\begin{proposition}
\begin{enumerate}
\item The generator $L^{c,n+1,n} $ defines a Markov process on
$\Omega_f\times\Omega_f$. For $\eta\geq\xi$ the path space measure
$P^{\eta,\xi}_{n+1,n}$ is a coupling of $P^\eta_{n+1}$ and $P^\xi_n$.
\item For $\eta\geq\xi $, 
$P^{\eta,\xi}_{n+1,n} (\eta_t\geq\xi_t,\ \forall t\geq 0)=1$.
\item For $f\in\mathcal{M}$, $\eta\in\Omega_f$, we have 
$S_{n+1} (t) f (\eta ) \geq S_n(t) f(\eta )$.
\end{enumerate}
\end{proposition}
{\bf Proof:} By construction any transition in the process with generator
$L^{c,n+1,n}$ preserves the order. To see point 3, 
for $f\in \mathcal{M}$
\begin{eqnarray}
S_{n+1} (t) f(\eta ) &=& \E^\eta_{(n+1)} f(\eta_t ) 
=\E^{\eta,\eta}_{n+1,n} f (\eta^1_t ) \geq
\E^{\eta,\eta}_{n+1,n} f(\eta^2_t ) \nonumber\\
&=& \E^\eta_n f(\eta_t) = S_n (t) f (\eta ).
\end{eqnarray}
\QED

\subsection{Construction of the process.}
We can now proceed to the definition of the sandpile process 
as a well-defined Markov process
with formal generator $L$. For functions $f\in \mathcal{M}$, on
$\Omega$
we define the semigroup associated to $L_n$ by
\begin{equation}\label{ho}
S_n(t) f (\eta ) = \lim_{\eta'\in \Omega_f,\eta'\uparrow\eta}S_n(t) f(\eta' ).
\end{equation}
In (\ref{ho}) the limit is taken along the increasing sequence
$\eta'_n$ defined by $\eta'_n (x) = \eta (x )$ for $|x|\leq n$, and
$\eta'_n (x)=1$ otherwise, but any other increasing sequence 
$\eta'_n\uparrow\eta$ gives the same limit.
Then, for $ f\in\mathcal{M}$ and $\eta\in\Omega$,
\begin{equation}\label{haha}
S_{n+1} (t) f (\eta ) \geq S_n(t) f (\eta ).
\end{equation}
Since $S_n(t)f (\eta)\leq \|f\|_\infty$
(where $\| f\|_{\infty} = \sup_{\eta\in\Omega} |f(\eta )|$), for 
$f\in\mathcal{M}$ and $\eta\in \Omega$ we can define
\begin{equation}\label{definitie}
S(t) f (\eta ) = \lim_{n\uparrow\infty} S_n (t) f (\eta ).
\end{equation}
Because each $S_n(t)$ is a semigroup on functions on $\Omega_f$, we have,
using (\ref{haha}), that
$S(t) 1=1$ and that $S(t)f\geq 0$ if $f\geq 0$. Moreover, for $f\in\mathcal{M}$, 

\begin{eqnarray}
S(t+s)f &=& \lim_{n\uparrow\infty} S_n (t+s) f
= \lim_{n\uparrow\infty} S_n(t) S_n (s) f\nonumber\\
&\leq& \lim_{n\uparrow\infty} S_n (t) (S(s)f)
= S(t)S(s) f,
\end{eqnarray}
and for all $m\in\N$
\begin{eqnarray}
S(t+s)f &=& \lim_{n\uparrow\infty} S_n(t)S_n(s)f
\geq \liminf_{n\uparrow\infty} S_m (t)( S_n (s)f)\nonumber\\
&\geq& S_m(t) (\liminf_{n\uparrow\infty} S_n(s) f)
=S_m (t) S(s)f,
\end{eqnarray}
so by taking the limit as $m\uparrow\infty$ we derive
the semigroup property for $S(t)$.
This semigroup 
defines the finite dimensional distributions of a Markov process;
by Kolmogorov's theorem, there is a unique Markov process 
$\{ \eta_t:t\geq 0\}$ such that $S(t)f(\eta ) = \E^\eta f(\eta_t )$. 
We call this process the {\it one-dimensional
sandpile process} (SP). We denote its path space measure
(i.e. a measure on $\Omega^{[0,\infty)}$) starting from $\eta$ by
$P_{SP}^{\eta}$.
\\[7mm]
{\bf Remarks.}

\begin{enumerate}
\item 
We cannot conclude that $P^\eta_{SP}$ has a version concentrating
on $D([0,\infty),\Omega)$ (the set of cadlag paths on $\Omega$) 
for all initial $\eta$ (cf. Corollary 5.1 below).  
This question is related to
the following one: For which $\eta\in\Omega$ do we have
\begin{equation}
\lim_{t\downarrow 0} S(t) f(\eta ) = f(\eta )
\end{equation}
for all local $f$? It is not even
clear that this limit exists. In contrast, this
is automatic for long
range exclusion: But in our case, because the values of
a configuration in a finite
volume are influenced by transitions outside this volume, the estimate
(2.3) of \cite{Ligg1} does not hold. 
\item 
We cannot expect that the formal generator $L$ is ``really" the
generator of the process. However we do prove in
Theorem 4.1 that for some class of ``decent"
configurations and local $f$
\begin{equation}\label{generator}
\lim_{t\downarrow 0} \frac{1}{t} (S(t)f(\eta)-f(\eta))= Lf(\eta ).
\end{equation}
This situation of pointwise convergence for a set of ``decent" configurations,
rather than uniform convergence, can be compared to the situation of
``weakly Gibbsian measures" where the potential is only convergent on a set
of good configurations, cf. \cite{MRV1},\cite{MRV2}. 
In the following section we will deal with
(\ref{generator}) and show that for short times (depending on $\eta$) and
for ``decent" configurations (a set of measure one for every ergodic probability
measure on $\Omega$, except $\delta_2$, the measure concentrating
on $\eta\equiv 2$) we actually have
\begin{equation}\label{stronger}
S(t) f(\eta ) = \sum_{n=0}^\infty \frac{t^n}{n!} L^n f(\eta),
\end{equation}
which is of course much stronger than (\ref{generator}) (but note that
we do have (\ref{stronger}) for the pure birth process).
\item 
In order to find an invariant measure for the sandpile process,
we cannot proceed in a standard
way. Since the process is not Feller (see
Corollary 5.1 below), a weak
limit point of $\frac{1}{T}\int_0^T \mu S(t) dt$ does not need to be invariant,
and $\int Lf d\mu =0$ does not imply that $\mu$ is invariant.
In Section 5 we shall show that the only invariant measure is
$\delta_2$ and that every initial measure converges to $\delta_2$
{\it in a finite time}. This explosive convergence to the invariant
measure is related
to the non-locality of the dynamics and to the fact that we continue to
add sand at rate one at each lattice site (cf. \cite{Carlson1}). 
\end{enumerate}

\section{Decent Configurations.}
As we mentioned in Section 2, the transformation $T_i$ is non-local. 
Nevertheless we still have
some kind of locality: Let $f:\Omega\rightarrow\R$ be a local function
depending on coordinates in $A\in \mathcal{S}$.
Put $A^- :=\min (A)$, $A^+:=\max (A)$, then we have
\begin{equation}
T_i f(\eta) - f(\eta) =0,\ \forall i\in \Z\setminus [A^--k^-(A^-,\eta),
A^++k^+(A^+,\eta)].
\end{equation}
This means that the effect of adding a sand particle is only felt in the
region bounded by the first one to the left and the first one
to the right of the
dependence set of $f$. Loosely speaking, $T_if$ has a dependence set which
is finite but dependent on the configuration.

We now formalize these statements. 
Let $\Omega_1$ denote the set of configurations
with an infinite number of ones to the right and to the left of the origin,
i.e.,
\begin{equation}
\Omega_1 := \{ \eta\in\Omega: \sum_{i<0} (2-\eta(i))=\sum_{i>0} (2-\eta (i))
=+\infty\}.
\end{equation}
For $\eta\in\Omega_1$ we order the sites $i\in\Z$ for which
$\eta (i) = 1$ as follows
\begin{equation}
\eta^{-1}(\{ 1\} ) := \{ X_i (\eta ):i\in\Z \}
\end{equation}
where $X_0 (\eta ):=\min \{ i\geq 0:\eta (i) =1 \}$ and $i<j$ implies
$X_i(\eta ) < X_j (\eta )$.
Then we define
\begin{eqnarray}
&I_0& = (X_{-1},X_0]\cap\Z\nonumber\\
&I_1& = (X_0,X_1]\cap\Z\nonumber\\
&I_{-1}& = (X_{-2},X_{-1}]\cap\Z,\ \mbox{ etc.}
\end{eqnarray}
In this way we  
view configuration $\eta$ as a sequence of intervals, where
each interval is bounded by two sites $i,j\in \Z, \ i<j$, for which
$\eta (i)= \eta (j ) = 1$. 
A function $f:\Omega_1\rightarrow\R$ will be called $N-${\it local} 
around $k \in \Z$ if
\begin{equation}
f(\eta ) = f\left(\eta(i), \ i\in\bigcup_{j=-N}^{N} I_{j+k} (\eta )\right).
\end{equation}
We say that $f$ is $N$-local if it is $N$-local around the origin.

\vspace{5mm}

{\bf Examples.}
\begin{enumerate}
\item $f(\eta) = \eta (i)$ is $0-$local around $i$, and every local function
is $N$-local for some $N$.
\item $f(\eta ) = |I_k (\eta )|$ is $k-$local for every $k$, but is
not local.
\item $f(\eta )= \sum_{x\in\Z} e^{-|x|}\eta (x)$, while continuous,
is not $N$-local for
any $N$.
\end{enumerate}

\begin{proposition}\label{boel}
Let $f:\Omega_1\rightarrow\R$ be $N-$local then
\begin{equation}
f(T_i\eta)-f(\eta) = 0 \mbox{ for all } i\in\Z\setminus
\bigcup_{j=-N-1}^{N+1} I_j (\eta ),
\end{equation}
and $Lf$ is ($N+1$)-local.
\end{proposition}
{\bf Proof:} If $i\in\cup_{j=-N-1}^{N+1} I_j (\eta )$, then, for
$f$ $N$-local, $f(T_i \eta )$ is $(N+1)$-local. Otherwise $T_i\eta
=\eta$. \QED

\begin{proposition}\label{grens}
Let $f:\Omega_1\rightarrow\R$ be $N-$local and bounded, then we have
\begin{equation}
|L^n f(\eta )|\leq
\left(|I_{-N-n} (\eta)| +\ldots +|I_{n+N} (\eta )|\right)^n 2^n \|f\|_\infty.
\end{equation}
\end{proposition}
{\bf Proof:} We put $N=0$ for the sake of simplicity, the case $N>0$ is treated
in the same way. We proceed by induction. For $n=1$
\begin{equation}
Lf(\eta ) = \sum_i [f(T_i\eta ) - f(\eta )],
\end{equation}
and, by Proposition \ref{boel}, this sum runs over 
$i\in I_{-1}(\eta )\cup I_0 (\eta )\cup I_1 (\eta )$,
$Lf$ is $1-$local, and 
\begin{equation}
|Lf(\eta ) | \leq (|I_{-1}(\eta)| + |I_0(\eta)| + |I_1(\eta)|)2\| f\|_\infty.
\end{equation}
Suppose now that for $n\geq 1$, $L^nf$ is $n-$local and
\begin{equation}
|L^n f(\eta ) | \leq (|I_{-n}(\eta )| + \ldots +|I_n(\eta )| )^n 
2^n \|f\|_\infty.
\end{equation}
Then we have
\begin{equation}\label{som}
L^{n+1} f(\eta ) = \sum_i [L^n f(T_i \eta)-L^n f(\eta )],
\end{equation}
where the sum runs over $i\in I_{-n-1}(\eta) \cup\ldots\cup I_{n+1}(\eta)$, 
since
$L^nf$ is $n-$local.
Therefore, for such an $i$,
\begin{eqnarray}
|L^n f (T_i\eta )| &\leq& 
\left[ |I_{-n} (T_i\eta )|+\ldots+|I_n (T_i\eta )|  \right]^n
2^n\| f\|_\infty\nonumber\\
&\leq& \left[ |I_{-n-1} (\eta )| +\ldots + |I_{n+1}  (\eta )|  \right]^n 2^n
\| f \|_\infty.
\end{eqnarray}
Thus
\begin{equation}
|L^{n+1} f(\eta ) | \leq 
\left[ |I_{-n-1}(\eta )|+\ldots + |I_{n+1} (\eta )|\right]^{n+1} 2^{n+1}
\| f \|_\infty.
\end{equation}
\QED
The following lemma is a direct consequence of Cauchy's formula for the
radius of convergence of a power series.
\begin{lemma}\label{cauchy}
Let $\{ a_n:n\geq 0 \}$ be a sequence of positive real numbers such that
$\limsup_{n\rightarrow\infty} a_n/n = a <\infty$. Then the series
$\sum_{n=0}^\infty t^n a_n^n/n!$ converges for $|t|<\frac{1}{a.e}$.
\end{lemma}
We now define the set of ``decent" configurations:
\begin{equation}\label{decent}
\Omega_{{\mbox {dec}}} := \{ \eta\in\Omega_1: a(\eta ):=
\limsup_{n\rightarrow\infty} \frac{|I_{-n}(\eta )|+\ldots+|I_n(\eta )|}{2n}
<\infty \}.
\end{equation}
Observe that for
$\eta\in\Omega_f$, if $|i|$ is large enough, $|I_i (\eta )|=1$,
hence $a(\eta ) =1$. Therefore 
$\Omega_f\subset\Omega_{{\mbox {dec}}}$.
More generally, if the density of ones
\begin{equation}
\rho (\eta ) =\lim_{n\rightarrow\infty} \frac{\sum_{i=-n}^n (2-\eta (i))}
{2n+1}
\end{equation}
exists, then $a(\eta ) = 1/\rho (\eta )$.
This implies that for every probability
measure $\mu$ ergodic under spatial translations and
for which $\mu (2-\eta (0)) >0$,
the set $\Omega_{{\mbox {dec}}}$ has 
$\mu$-measure 1. It is also clear that $\Omega_{{\mbox {dec}}}$ is 
a translation invariant set in the tail field, 
i.e., if $\eta\in\Omega_{{\mbox {dec}}}$ and $\zeta$ differs
from $\eta$ in a finite number of lattice sites, then  
$\zeta\in\Omega_{{\mbox {dec}}}$.

\begin{theorem}\label{serie}
Let $f$ be $N-$local for some $N\in\N$ and $\eta\in\Omega_{{\mbox {dec}}}$, then
for $t<\frac{1}{4ea(\eta)}$ the series 
$\sum_{n=0}^\infty \frac{t^n L^n f(\eta )}{n!}$ converges
absolutely and equals $S(t)f(\eta)$. In particular, $L$ is the
``pointwise generator" of the semigroup, i.e.
\begin{equation}
\lim_{t\downarrow 0} \frac{S(t)f(\eta)-f(\eta)}{t} = Lf (\eta ).
\end{equation}
\end{theorem}
{\bf Proof:} The absolute convergence of the series follows from
Proposition 4.2, Lemma 4.1 and (\ref{decent}).
By definition, for an $N$-local $f\in\mathcal{M}$ 
\begin{equation}
S(t) f(\eta )= \lim_{n\uparrow\infty} \lim_{\eta'\in\Omega_f,\eta'\uparrow\eta}
S_n(t) f(\eta' ).
\end{equation}
For $\eta'\in\Omega_f$, $a(\eta') =1$ and 
\begin{equation}\label{dada}
S_n (t) f (\eta' ) = \sum_{k=0}^\infty t^k/k! (L_n)^kf (\eta' ).
\end{equation}
Moreover, for every $ k\in\N,\  \eta\in\Omega$:
\begin{equation}
\lim_{n\uparrow\infty}\lim_{\eta'\in\Omega_f,\eta'\uparrow\eta}
(L_n)^k f(\eta' )=L^k f(\eta ).
\end{equation}
We can bring in the
limits $\lim_{n\uparrow\infty}\lim_{\eta'\uparrow\eta}$ into the sum
of (\ref{dada}) by the dominated convergence theorem 
together with Lemma 4.1 and the inequalities
\begin{eqnarray}
|L^k_n f(\eta' )| &\leq & \left(
|I_{-N-k} (\eta' )| +\ldots + |I_{N+k} (\eta' ) | \right)^k 2^k
\| f\|_{\infty}\nonumber\\
&\leq & \left( |I_{-N-k} (\eta ) + \ldots +|I_{N+k} (\eta )|
\right)^k 2^k \| f\|_{\infty},
\end{eqnarray}
for $\eta'\leq\eta,\ \eta'\in\Omega_f$.\QED

{\bf Remark.}
Theorem \ref{serie} shows that the infinite volume sandpile
process is indeed a natural extension of the finite volume
standard abelian sandpile model. More precisely, if $P_n$ denotes
the transition probability matrix of the discrete-time sandpile
model in the finite volume
$\Lambda_n=[-n,n]\cap\Z$,
then we have for $\eta\in\Omega_{\mbox{dec}}$ and $t$
small enough (depending on $\eta$),
\begin{equation}
\lim_{n\uparrow\infty} P_n^{\lfloor nt \rfloor} f(\eta^n ) = S(t) f(\eta ),
\end{equation}
where $\eta^n$ denotes the restriction of $\eta$ to $\Lambda_n$. 
\section{Invariant measure for the sandpile process.}
In this section we show that after a finite time, the process reaches
$\delta_2$, which is therefore its only invariant measure.
Nevertheless we first prove Proposition 5.1 to illustrate
some ``pathological" behavior of the process (cf. Corollary 5.1).

\begin{proposition}\label{deltatwee}
$\delta_2$ is invariant for the sandpile process.
\end{proposition}
{\bf Proof:} 
Denote by $\hat{2}$ the configuration $\eta\equiv 2$. 
Notice that this is not a decent configuration. Thus we have to
prove that
\begin{equation}\label{two}
S(t)f(\hat{2})=f(\hat{2}).
\end{equation}
Without loss of generality we put $f(\eta ) = \eta (0)-1$. 
Clearly, $S(t) f(\hat{2} ) \leq f(\hat{2})$, so we are left to show the
opposite inequality. By definition (\ref{definitie}),

\begin{eqnarray}
S(t)f(\hat{2}) &=& \lim_{n\rightarrow\infty}
\lim_{\eta'\in\Omega_f,\eta'\uparrow\eta} S_n(t)f(\eta')
\nonumber\\
&\geq & \lim_{\eta'\in\Omega_f,\eta'\uparrow\eta} S_a (t) f(\eta' )
\nonumber\\
&=& \lim_{n\uparrow\infty} P^{A_n} (0\in A(t)).
\end{eqnarray}
Here $P^{A_n}$ denotes the Markov measure on $\mathcal{S}^{[0,\infty)}$
associated to the process with generator $\hat{L_a}$ (only avalanches),
starting from the set $A_n:= [-n,n]\cap\Z$. The set $A(t)$ will always be 
of the form $I(t)\setminus\{ x\}$
where $I(t)$ is an interval which grows in time.
The probability that $\{ x \} = \{ 0 \}$ will thus be of the order
$(|A_n| + f(t))^{-1}$ where $f(t)$ is some increasing function of $t$.
Let us make this precise. Denote by $\tau_1,\tau_2,\ldots,\tau_n$ the successive
event times of the Markov chain $\{ A(t) : t\geq 0 \}$. Then we have

\begin{eqnarray}
P^{A_n} (0\not\in A(\tau_1) ) &=& P^{A_n} (x\not\in A(\tau_1) ), 
\ \forall x\in A_n\nonumber\\
&=& \frac{1}{|A_n|},\nonumber\\
P^{A_n} (0\not\in A(\tau_2) ) &=& \sum_{B\in \mathcal{S}}
\frac{1}{|B|} P^{A_n} (A(\tau_1) = B)= \frac {1}{|A_n|+1},
\end{eqnarray}
and analogously

\begin{equation}
P^{A_n} (0\not\in A(\tau_n) ) = \frac{1}{|A_n| + n-1}.
\end{equation}
Hence we conclude that for all $t\geq 0$

\begin{equation}
P^{A_n} (0\not\in A (t)) 
\leq \frac{1}{|A_n|} = \frac{1}{2n+1}.
\end{equation}
Therefore

\begin{equation}
\limsup_{n\uparrow\infty} P^{A_n} (0\not\in A (t) ) = 0
\end{equation}
which proves the claim. \QED

\begin{corollary}
The sandpile process is not Feller, and for some initial
configurations the process has no right-continuous version.
\end{corollary}
{\bf Proof:} As in the preceding proof, 
we have that for all $\eta\in\Omega$
such that $\eta^{-1}(\{1\} ) \in \mathcal{S}$ (finite number of ones)
\begin{equation}
P^{\eta}_{SP} (\eta_t (x) \not= 2) = 0, \ \forall t>0.
\end{equation}
Therefore there is no right-continuous version of the process starting
from these configurations (cf. Remark 1 of Section 3.4).\QED
The following theorem will be proven at the end of this section.

\begin{theorem}
In the sandpile process there exists a finite $T>0$ such 
that for all $t\geq T$, $\eta \in\Omega$,
\begin{equation}
P_{SP}^\eta [\eta_t (0) = 1 ] = 0.
\end{equation}
\end{theorem}

 This result seems intuitive: Since
we continue adding sand at each lattice site with rate $1$, after a time
of order unity, every site will have received an extra sand particle.
Once the configuration is full, it remains full. Nevertheless, this intuitive
picture is somewhat dangerous, since we are working in infinite volume,
and it could be possible that sand particles ``disappear" at infinity.
Indeed, the same statement is not true 
anymore in higher dimensions: Although we still
expect convergence to the invariant measure in a finite
time, this invariant measure will not be the Dirac measure concentrating on
the ``full configuration". 

\subsection{The finite volume sandpile process.}
For the finite volume sandpile process in $[-n,n]$, we only add
sand particles in $[-n,n]$,
and at the boundaries we put fixed ones (i.e., sand falls off the
boundary and disappears when the boundary site topples). 
More formally, for $\eta\in\Omega$, we
denote by $\eta^n$ the configuration which is equal to
$\eta$ in $[-n,n]$ and is identically
$1$ outside that interval. Define
\begin{equation}
T_{n,i} \eta := \left[ T_i\eta^n \right]^n,
\end{equation}
\begin{equation}
\Omega_{f,n} := \{ \eta\in\Omega: \eta (j)=1,\ \forall j\in\Z\setminus[-n,n] \}.
\end{equation}
Note that $\Omega_f =\cup_{n} \Omega_{f,n}$.
The finite volume sandpile process (FVSP) is the pure jump process on 
$\Omega_{f,n}$
with generator
\begin{equation}\label{geny}
\mathcal{L}_nf(\eta ) = \sum_{i=-n}^n [f(T_{n,i} \eta ) - f (\eta ) ]
\end{equation}

\subsection{Poisson representation.}
Let $\{ N_t:t\geq 0 \}:=\{ N^i_t, i\in\Z, t\geq 0 \}$ be a collection of
independent rate one Poisson processes. For $\Xi\in \N_0^{\Z}$
we define $\theta_n (\Xi )$ to be the configuration obtained from
``toppling" in $[-n,n]$ and fixing ones outside, i.e.,
\begin{equation}\label{teita}
\theta_n (\Xi ) = \left[ \prod_{i=-n}^n T^{\Xi(i)-1}_{n,i} \right]\hat{1}
\end{equation}
where $\hat{1}$ denotes the configuration $\eta\equiv 1$, and 
$T^{\Xi(i)-1}_{n,i}$ means $\Xi(i)-1$ iterations of $T_{n,i}$.
\begin{proposition}
The process $\{ \theta_n (N_t +\eta ) :t\geq 0\}$ is a version of
the FVSP starting from $\eta\in\Omega_{f,n}$ ($N_t+\eta$ means coordinatewise
addition).
\end{proposition}
{\bf Proof:} 
Notice that the transformations $T_{n,i}$ commute for different $i$. Therefore,
first adding all the particles and then letting the system topple is equivalent
to adding the particles one by one and letting the system topple each time.
More formally, using (\ref{teita}), 
\begin{equation} 
\left(\frac{d}{dt} \E f(\theta_n (N_t + \eta))
\right)_{t=0}=
\mathcal{L}_nf(\eta )
\end{equation} 
for $\eta\in\Omega_{f,n}$.\QED

\subsection{Invariant measure for the FVSP.}
Since the FVSP is a pure jump process on the finite set $\Omega_{f,n}$, it
converges exponentially fast to its unique stationary measure $\mu_n$. 
\begin{equation}
\mu_n = \frac{1}{|\mathcal{R}_n|} \sum_{\eta\in\mathcal{R}_n} \delta_\eta
\end{equation}
is the uniform measure on the set $\mathcal{R}_n$ of recurrent configurations: 
\begin{equation}
\mathcal{R}_n =  \{ \eta\in\Omega_{f,n}:\sum_{i=-n}^n (2-\eta (i))\leq 1\},
\end{equation}
i.e., those configurations in $\Omega_{f,n}$ containing at most one site
$i\in [-n,n]$ such that $\eta (i) = 1$. For all
$i\in [-n,n]\cap\Z$, $\mu_nT_{n,i} = \mu_n$. More precisely for every pair
$(\eta,\zeta)\in\mathcal{R}_n\times\mathcal{R}_n$, $\eta\not=\zeta$ 
there exists a unique
lattice site $i\in [-n,n]$ such that $T_{n,i} \eta = \zeta$. Therefore,
\begin{eqnarray}
\int (\mathcal{L}_nf)g d\mu_n &=& \sum_{\eta\in\mathcal{R}_n}
\frac{1}{|\mathcal{R}_n|}\sum_{i = -n}^n
\left[ f(T_{n,i}\eta ) - f(\eta ) \right] g(\eta )\nonumber\\
&=&\frac{1}{|\mathcal{R}_n|}\sum_{\eta,\zeta\in\mathcal{R}_n,\eta\not=\zeta} 
\left[ f(\zeta ) - f(\eta ) \right]
g(\eta )\nonumber\\
&=& \int f(\mathcal{L}_ng) d\mu_n.
\end{eqnarray}
so that $\mu_n$ is actually reversible for $\mathcal{L}_n$.

\subsection{Convergence to equilibrium.}
The key inequality in the proof of Theorem 5.1 is: There exists
 $T>0$ such that for all $t\geq T$
\begin{equation}
P(\theta_n (N_t + \eta ) (0) = 1 ) \leq Ce^{-n\lambda} +C'/n.
\end{equation}
Indeed, after a finite time
$T$ in the Poisson process we have for large $n$
\begin{equation}
\sum_{i=-n}^n N^i_T\approx 2nT.
\end{equation}
For $T$ large enough this represents a serious ``excess" at most lattice
sites in $[-n,n]$. If we let this configuration topple, then we will
obtain a recurrent configuration, as we show in the
following lemma.

\begin{lemma}\label{geprul}
Let $\eta\in\N_0^{\Z}$ be such that 
\begin{equation}\label{overdreven}
\sum_{i = -\lfloor n/2 \rfloor }^{\lfloor n/2 \rfloor} \eta (i) \geq 12n,
\end{equation}
then $\theta_n (\eta ) \in \mathcal{R}_n$.
\end{lemma}
{\bf Proof:} First of all note that we put the ``excess" of particles sufficiently
far from $-n$ and $n$, so that not too much sand
will be lost at the boundaries. We call a site $i\in [-n,n]\cap\Z$ an
{\it excess site} if $\eta (i) >2$. For the sake of simplicity we 
suppose $n$ even. Put
\begin{eqnarray}
m &=& \min \{ i\in [-n/2,n/2]\cap\Z :\eta (i) >2 \}\nonumber\\
M &=& \max \{ i\in [-n/2,n/2]\cap\Z :\eta (i) >2 \}.
\end{eqnarray}
Because of the monotonicity property
\begin{equation}
\eta,\zeta\in{\N_0}^{\Z},\ \zeta\geq\eta,\mbox{ and } 
\theta_n (\eta )\in\mathcal{R}_n
\mbox{ imply } \theta_n (\zeta ) \in \mathcal{R}_n,
\end{equation}
we can replace our configuration $\eta$ by $\eta'$ defined by:
$\eta'(i)=\eta (i)$ for $i\in[-n/2,n/2]\cap\Z$ and $\eta' (i)=1$ otherwise.
We now let the system topple in the following way: All sites in $(m,M)$
topple and the sand arriving at $m$ or $M$ will be stocked at these sites
(i.e. sites $m$ and $M$ do not topple). We then end up with a
configuration $\eta''$ with two excess sites at $m$ and $M$ and since
no sand is lost in this toppling:
\begin{equation}
\eta''(m) + \eta'' (M) \geq 10n.
\end{equation}
We change again the configuration only retaining the site with the
maximum excess, say $M$, and replacing all the rest by 1's, this configuration
is denoted by $\eta'''$. The excess of $\eta'''$ at $M$ is not less than
$5n$. So we are left to prove that a configuration with one excess site of
this type topples to a recurrent configuration. The case $M=0$ is trivial.
Suppose $M>0$. 
After $n-M$ topplings
we get a configuration $\eta^{iv}$ with an excess site at $M$
(an excess not less than $3n$), $\eta^{iv} (i) =2$ for all $i\in [2M-n,n]$,
$\eta^{iv} (i) =1$mbox for all $i\in \Z\setminus [2M-n,n]$. From that moment on,
sand particles are lost at the boundary by further toppling. It is easy to see
that the worst situation is when $M=n/2$. In that case, define the configuration 
$\eta^{v}$ by
$\eta^{v} (i)=2$ for $i\in [0,n]$ and $\eta^v (i) =1$ otherwise; since
the excess at site $i=n/2$ in $\eta^{iv}$ 
is not less than $3n$, it suffices to
show that
$(T_{n,n/2})^{3n} (\eta^v ) \in \mathcal{R}_n$. 

By adding sand particles
at site $n/2$ in $\eta^v$, we shall always obtain a configuration in the
set $\Gamma_n=\cup_{-n\leq a<0} \Gamma_{n,a}$, with
\begin{equation}
\Gamma_{n,a}=\{ \eta\in\Omega_{f,n}: \ 
\eta (i) =1 \ \forall i<a, \mbox{ and } \sum_{i=a}^n (2-\eta (i))\leq 1 \}.
\end{equation}
That is, for $\eta\in\Gamma_{n,a}$,  
there exists a ``leftmost" site $a<0$ for which $\eta (a)=2$ and
the interval $[a,n]$ contains at most one site $i$ for which $\eta (i)=1$.
When $a$ is the position of the leftmost two in the configuration $\eta$, then
by adding three sand particles at site $n/2$ and letting the system topple,
the position of the leftmost two becomes $a-1$. I.e.,
for $a<0$, $\eta\in\Gamma_{n,a}$ implies $(T_{n,n/2})^3\in \Gamma_{n,a-1}$.
By this last property, since $\Gamma_{n,-n}=\mathcal{R}_n$, 
$(T_{n,n/2})^{3n}\eta^{v}\in\mathcal{R}_n$.
\QED

The following lemma is a standard large deviation estimate for sums
of independent random variables. 

\begin{lemma}\label{largedev}
For every strictly positive 
$t$ and $ \epsilon $, there exists  a constant $\lambda (\epsilon ) >0$
such that
\begin{equation}
P\left(\sum_{x=- \lfloor n/2 \rfloor}^{\lfloor n/2 \rfloor} 
N^x_t \leq n(t-\epsilon)\right)\leq 
\exp{(-n\lambda(\epsilon))}.
\end{equation}
\end{lemma}
{\bf Proof of Theorem 5.1:} We show that the statement is true for $T=14$. 
By a coupling argument as in Section 3, we get that
the FVSP $\theta_n (N_t+\eta )$ is stochastically dominated by the sandpile
process $\eta_t$ (since we put fixed ones at the boundaries of the interval
$[-n,n]$ in the FVSP). Put $0<\delta <1$. By Lemmas \ref{largedev}
and \ref{geprul} we know that for $t>14$, there exists
$\lambda>0$ such that
\begin{equation}
P (\theta_n (N_{t-\delta}) \not\in \mathcal{R}_n)\leq 
C\exp{(-\lambda n)}.
\end{equation}
Hence we obtain the inequalities
\begin{eqnarray}\label{afsch}
P_{SP}^\eta (\eta_t (0) = 1 ) &\leq & P (\theta_n (N_t + \eta ) (0)
=1)\nonumber\\
&\leq & P [ \theta_n (N_t +\eta ) (0) =1|\theta_n (N_{t-\delta}+\eta )\in
\mathcal{R}_n]\nonumber\\ 
&+& C\exp{(-\lambda n)}.
\end{eqnarray}
Conditioned on $\{\zeta:=\theta(N_{t-\delta}+\eta )\in\mathcal{R}_n\}$,  
the distribution of $\theta_n (N_t +\eta )$ is  uniform
on a set of the form $\mathcal{R}_n\setminus\{ \zeta' \}$ if there has been
at least one Poisson event in $(t-\delta,t)$. 
Therefore, since there is only one $\eta\in\mathcal{R}_n$
with $\eta (0) =1$,
\begin{eqnarray}\label{gruwel}
&P& \left[ \theta_n (N_t+\eta ) (0) =1|\theta_n (N_{t-\delta}+\eta)
i\in\mathcal{R}_n \right]\nonumber\\
&\leq & \frac{1}{2n+1} + (1-e^{-\delta})^{2n+1}.
\end{eqnarray}
Combining (\ref{gruwel}) with (\ref{afsch}) we arrive at
\begin{equation}\label{dawaset}
P^\eta_{SP} [\eta_t (0) =1 ] \leq \frac{1}{2n+1} + (1-e^{-\delta})^{2n+1}
+ Ce^{-n\lambda},
\end{equation}
which concludes the proof of the theorem.\QED
\\[7mm]
{\bf Acknowledgements:} This work was financially supported by
the agreement ``Tournesol, no. 97.049". We thank Herv\'{e} Guiol for useful
comments.

\end{document}